\documentclass[11pt]{article}

\usepackage{latexsym,mathrsfs}
\usepackage{amsmath,amssymb}
\usepackage{amsthm,enumerate,verbatim}
\usepackage{amsfonts}
\usepackage{graphicx}
\usepackage{algorithm}
\usepackage{algorithmic}
\usepackage{url}

\newtheorem{lemma}{Lemma}

\newtheorem{theorem}{Theorem}
\newtheorem{corollary}{Corollary}
\newtheorem{proposition}{Proposition}
\newtheorem{definition}{Definition}
\newtheorem*{definition*}{Definition}
\newtheorem{example}{Example}

\newtheorem{remark}{Remark}
\newtheorem*{problem*}{Problem}

\numberwithin{equation}{section}
\numberwithin{table}{section}
\numberwithin{figure}{section}

\DeclareMathOperator{\argmin}{argmin}

\DeclareMathOperator{\tr}{tr}

\newcommand {\mat}  [1] {\left[\begin{array}{#1}}
\newcommand {\rix}      {\end{array}\right]}

\def\real{\mathop{\mathrm{Re}}}

\newcommand{\eproof}{\space
    {\ \vbox{\hrule\hbox{\vrule height1.3ex\hskip0.8ex\vrule}\hrule}}\par}

\def \R{{\mathbb R}}
\def \C{{\mathbb C}}

\title{Computing nearest stable matrix pairs}

\author{Nicolas Gillis\thanks{Department of Mathematics and Operational Research, University of Mons, Rue de Houdain 9, 7000 Mons, Belgium. Email: \{nicolas.gillis, punit.sharma\}@umons.ac.be. NG and PS acknowledge the support of the ERC (starting grant n$^\text{o}$ 679515). NG also acknowledges the support of the F.R.S.-FNRS (incentive grant for scientific research n$^\text{o}$ F.4501.16). }
\and Volker Mehrmann\thanks{Institut f${\rm \ddot{u}}$r Mathematik, MA 4-5 TU Berlin, Str.\@ d.\@ 17.\@ Juni 136, D-10623 Berlin, Germany.
Email: mehrmann@math.tu-berlin.de.
VM acknowledges support by the DFG Collaborative Research Center 910 {\it Control of self-organizing nonlinear systems: Theoretical methods and concepts of application.}} 
\and Punit Sharma$^*$}

\begin{document}

\maketitle

\begin{abstract}
In this paper, we study the nearest stable matrix pair problem: given a square matrix pair $(E,A)$,
minimize the Frobenius norm of $(\Delta_E,\Delta_A)$ such that $(E+\Delta_E,A+\Delta_A)$ is a stable matrix pair.
We propose a reformulation of the problem with a simpler feasible set  by
introducing dissipative Hamiltonian (DH) matrix pairs: A matrix pair $(E,A)$ is DH if $A=(J-R)Q$ with skew-symmetric $J$, positive semidefinite $R$, and an invertible $Q$ such that $Q^TE$ is positive semidefinite. This reformulation has a convex feasible domain onto which it is easy to project. This allows us to employ a fast gradient method to obtain a nearby stable approximation of a given matrix pair.
\end{abstract}


\section{Introduction}
We study the stability and instability under perturbations for systems of linear time-invariant differential-algebraic equations (DAEs) of the form
\begin{equation}\label{eq:lindae}
E\dot{x}=Ax+f,
\end{equation}
on the unbounded interval ${\mathbb I}=[t_0,\infty)$,
where $E,A \in \mathbb R^{n,n}$ and $f\in C^0({\mathbb I},\mathbb R^n)$ is sufficiently smooth.
Here $\mathbb R^{n,n}$ denotes the real $n\times n$ matrices and $C^k(D,\mathbb R^n)$ denotes the $k$-times differentiable functions from a set $D$ to $\mathbb R^{n}$.
Systems of the form~(\ref{eq:lindae}) arise from linearization around stationary solutions
of initial value problems for general implicit systems of DAEs
\begin{equation}\label{DAE}
F(t,x,\dot x)=0,
\end{equation}
with an initial condition
\begin{equation}\label{DAEIC}
x(t_0)=x_0.
\end{equation}
Here we assume that
$F\in C^0({\mathbb I}\times{\mathbb D}_x\times
{\mathbb D}_{\dot x},{\mathbb R}^n)$ is sufficiently smooth and
${\mathbb D}_x,{\mathbb D}_{\dot x}\subseteq{\mathbb R}^n$ are open sets, and
following  \cite{KunM06}, we use the following solution and stability concepts.
\begin{definition}\label{def_1}
Consider system (\ref{DAE}) with sufficiently smooth $F$.
\begin{itemize}
\item [i)]
A function $x:{\mathbb I} \rightarrow\mathbb R^n$
is called a {\em solution} of (\ref{DAE}) if
$x\in C^1({\mathbb I},\mathbb R^n)$
and $x$ satisfies (\ref {DAE}) pointwise.
\item [ii)] It  is called a {\em solution of the initial value problem}
(\ref{DAE})--(\ref{DAEIC})
if $x$ is a solution of (\ref {DAE}) and satisfies (\ref{DAEIC}).
\item [iii)] An initial condition (\ref {DAEIC}) is called {\em consistent}
if the corresponding initial value problem has at least one solution.
\item [iv)]
System (\ref{DAE}) is called regular if for every consistent initial condition $x_0$ it has a unique solution.
\end{itemize}
\end{definition}
\begin{definition}\label{stabtraj}
Suppose that the DAE~(\ref{DAE}) is regular and (\ref{DAEIC}) is a consistent initial condition.
A solution~$x:t\mapsto x(t;t_0,x_0)$ of (\ref{DAE})--(\ref{DAEIC}) is called
\begin{enumerate}
\item [i)]
\emph{stable} if for every $\varepsilon >0$ there exists $\delta>0$ such that
\begin{enumerate}
\item
the initial value problem (\ref{DAE}) with initial condition
$x(t_0)=\hat x_0$ is solvable on~${\mathbb I}$ for all consistent $\hat x_0\in {\mathbb R}^n$
with $\| \hat x_0 -x_0 \|<\delta$;
\item
the solution $x(t;t_0,\hat x_0)$ satisfies
$\| x(t;t_0,\hat x_0) -x(t;t_0,x_0) \|<\varepsilon$ on~${\mathbb I}$.
\end{enumerate}
\item [ii)]
\emph{asymptotically stable} if
it is stable and there exists $\varrho>0$ such that
\begin{enumerate}
\item
the initial value problem (\ref{DAE}) with initial condition
$x(t_0)=\hat x_0$ is solvable on~${\mathbb I}$ for all consistent $\hat x_0\in {\mathbb R}^n$
with $\| \hat x_0 -x_0 \|<\varrho$;
\item
the solution $x(t;t_0,\hat x_0)$ satisfies
$\lim_{t\to \infty}\| x(t;t_0,\hat x_0) -x(t;t_0,x_0) \|=0$.
\end{enumerate}
\end{enumerate}
\end{definition}
To characterize regularity and stability for the linear constant coefficient case we introduce
the following notation.
A square matrix pair $(E,A)$ with $E,A \in \mathbb R^{n,n}$ is called \emph{regular} if the matrix pencil $z E-A$ is regular, i.~e. if
$\operatorname{det}(\lambda E-A)\neq 0$ for some $\lambda \in \mathbb C$, otherwise it is called \emph{singular}.
For a regular matrix pair $(E,A)$, the roots of the polynomial $\operatorname{det}(z E-A)$ are called \emph{finite eigenvalues} of the pencil
$zE-A$ or of the pair $(E,A)$, i.~e. $\lambda\in \mathbb C$ is a finite eigenvalue of the pencil
$zE-A$ if there exists a vector $x\in \mathbb C^n\setminus \{ 0\}$ such that
$(\lambda E-A)x=0$, and $x$ is called an \emph{eigenvector} of $zE-A$ corresponding to the eigenvalue $\lambda$.
A regular pencil $zE-A$ has \emph{$\infty$ as an eigenvalue} if $E$ is singular.

A regular real matrix pair $(E,A)$ (with $E,A\in \R^{n,n}$) can be transformed to \emph{Weierstra\ss\ canonical form} \cite{Gan59a}, i.~e. there exist nonsingular matrices $W, T \in \C^{n,n}$ such that
\[
E=W\mat{cc}I_q& 0\\0&N\rix T \quad \text{and}\quad A=W \mat{cc}J &0\\0&I_{n-q}\rix T,
\]
where $J \in \C^{q,q}$ is a matrix in \emph{Jordan canonical form} associated with the $q$ finite eigenvalues of
the pencil $z E-A$ and  $N \in \C^{n-q,n-q}$ is a nilpotent matrix in Jordan canonical form  corresponding
to  $n-q$ times the  eigenvalue $\infty$. If $q < n$ and $N$ has degree of nilpotency $\nu \in \{1,2,\ldots\}$, i.~e.
$N^{\nu}=0$ and $N^i \neq 0$ for $i=1,\ldots,\nu-1$, then $\nu$ is called the \emph{index of the pair} $(E,A)$. If $E$
is nonsingular, then by convention the index is $\nu=0$.
A pencil $zE-A$ is of index at most one  if it is regular with exactly $r:=\text{rank}(E)$ finite eigenvalues, see e.~g. \cite{Meh91,Var95}. In this case the $n-r$ copies of the eigenvalue $\infty$ are semisimple (non-defective).

Then we have the following well-known results, see e.~g.~\cite{ByeN93,DuLM13,KunM07}.
\begin{proposition}\label{prop1}
Consider the initial value problem~\eqref{eq:lindae}.
\begin{itemize}
\item [i)] The linear DAE~\eqref{eq:lindae} is regular
if and only if the matrix pair $(E,A)$ is regular~\cite{KunM06}.
\item [ii)] If the pair $(E,A)$
is regular and of index at most one, then the initial value problem~\eqref{eq:lindae}, \eqref{DAEIC} with a consistent initial value $x_0$ is stable if all the finite eigenvalues of $zE-A$ are in the closed left half of the complex plane and those on
the imaginary axis are semisimple.
\item [iii)] If the pair $(E,A)$
is regular and of index at most one then the initial value problem~\eqref{eq:lindae}, \eqref{DAEIC} with a consistent initial value $x_0$ is \emph{asymptotically stable} if  all the finite eigenvalues of $zE-A$ are in the open left half of the complex plane.
\end{itemize}
\end{proposition}
The literature on (asymptotic) stability of constant coefficient DAEs is very ambiguous, see e.~g. \cite{BoyGFB94,ByeN93,Var95}, and the review in \cite{DuLM13}. This ambiguity arises from the fact that some authors consider only the finite eigenvalues in the stability analysis and allow the index of the pencil $(E,A)$ to be arbitrary, others consider regular high index pencils $zE-A$ as unstable by considering $\infty$ to be on the imaginary axis. 

In this paper, we are interested in the problem of finding minimal perturbations to the system matrices $E+\Delta E$ and $A+\Delta A$ of an unstable DAE~\eqref{eq:lindae} that puts the system on the boundary of the stability region. This problem is closely related to the classical state-feedback-stabilization problem in descriptor control problems with $f=Bu$ for some control function $u$, see e.~g. \cite{Ben11,Dai89,Meh91}, where one chooses feedbacks $u=Kx$ so that the closed-loop descriptor system
\[
E \dot x =(A+BK)x
\]
is asymptotically stable. However, then  the perturbations are restricted only in $A$ while in derivative feedback
one also allows feedbacks $E+BG$, see \cite{BunMN92,BunMN94}.
Instead of the feedback stabilization problem, we extend the nearest stable matrix problem of  \cite{GilS16,OrbNV13} to matrix pencils and we study the following problem.
\begin{problem*}{\rm
For a given  pair $(E,A) \in \R^{n,n}\times \R^{n,n}$
find the nearest asymptotically stable matrix pair $(M,X)$. More precisely,  let $\mathbb S$ be the set of
matrix pairs $(M,X)\in \R^{n,n}\times \R^{n,n}$ that are regular, of index at most one, and have all finite eigenvalues in the open left half plane,
then we wish to compute
\begin{equation*}
\inf_{(M,X) \in \mathbb S} \{{\|E-M\|}_F^2+{\|A-X\|}_F^2\},
\end{equation*}
where ${\|\cdot\|}_F$ stands for the Frobenius norm of a matrix.
}
\end{problem*}
In the following we use \emph{nearest stable matrix pair problem} to refer to the above problem.
This problem is the complementary problem to the distance to instability for matrix pairs, see~\cite{ByeN93} for complex pairs, and~\cite{DuLM13} for a survey on this problem.
Since we require a stable pair to be regular it is also complementary to the distance to the nearest singular pencil, which is a long-standing open problem \cite{ByeHM98,GugLM16,MehMW15}.

For ordinary differential equations, i.~e. when $E=I_n$,
stability solely depends on the eigenvalues of the matrix $A$. Recently, in \cite{GilS16} a new algorithm was proposed to compute a nearby stable perturbation to a given unstable matrix by reformulating this highly nonconvex problem into an equivalent optimization problem with a convex feasible region. This reformulation
uses the fact that the set of stable matrices can be characterized as the set of \emph{dissipative Hamiltonian (DH) matrices}. We recall from \cite{GilS16} that a  matrix $A \in \R^{n,n}$ is called a DH matrix if $A$ can be factored as $A=(J-R)Q$ for some $J,R,Q \in \R^{n,n}$ such that
$J^T=-J$, $R=R\succeq 0$, i.~e. is positive semidefinite, and $Q=Q^T \succ 0$, i.~e. is positive definite. In this way $A$ is stable if and only if $A$ is a DH matrix, see also \cite{BeaMX15_ppt,MehMS16,MehMS17}.

In this paper we extend the ideas of \cite{GilS16} from the matrix case to the pencil case and discuss the problem of computing the nearest stable matrix pair. This is a challenging problem, since the set $\mathbb S$ of all stable matrix pairs is neither open nor closed with respect to the usual norm topology as the following example demonstrates.
\begin{example}\label{ex:ex1}{\rm
Consider the matrix pair
\[
(E,A)=\left (\mat{ccc}1 &0&0\\0&0&0 \\ 0&0&0\rix,\ \mat{ccc}-1 &0&2\\0&1&0\\0&0&1\rix\right ).
\]
It is easy to check that $(E,A)$ is regular, of index one, and asymptotically stable with the only finite eigenvalue $\lambda_1=-1$, and thus
$(E,A) \in \mathbb S$. Consider the perturbations $(\Delta_E,\Delta_A)$,
where
\[
\Delta_E=\mat{ccc}\epsilon_1&0&0\\ 0&\epsilon_2 & \epsilon_3\\0&0&0\rix \quad \text{and}\quad
\Delta_A=\mat{ccc} \delta &0&0\\0&0&0\\0&0&0\rix,
\]
and the perturbed pair $(E+\Delta_E,A+\Delta_A)$. If we let $\delta=\epsilon_1=\epsilon_2=0$ and $\epsilon_3 >0$,
then the perturbed pair is still regular and all finite eigenvalues are in open left half plane, but it is of index two. If we let $\delta=\epsilon_1=\epsilon_3=0$ and $\epsilon_2 > 0$, then the perturbed pair is regular and index one, but has
two finite eigenvalues $\lambda_1=-1$ and $\lambda_2=1/\epsilon_2 >0$, and thus is not stable. This shows that
$\mathbb S$ is not open.
If we let $\epsilon_1=-\delta$, $\epsilon_2=\epsilon_3=0$, then as $\delta \rightarrow 1$ the perturbed
pair becomes non-regular. This shows that $\mathbb S$ is not closed.
}
\end{example}
It is also clear that the set $\mathbb S$ is highly nonconvex, see~\cite{OrbNV13}, and thus it is very difficult to
find a globally optimal solution of the nearest stable matrix pair problem. To address this challenge, we follow the strategy suggested in~\cite{GilS16} for matrices and reformulate the problem of computing the nearest stable matrix pair by extending the concept of dissipative Hamiltonian matrices to \emph{dissipative Hamiltonian matrix pairs (DH pairs)}.

The paper is organized as follows. In Section 2, we introduce and study DH matrix pairs. We provide several theoretical results and
characterize the set of asymptotically stable matrix pairs in terms of DH matrix pairs. A reformulation
of the nearest stable matrix pair problem using the DH characterization is obtained in
Section 3. In Section 4, we propose a fast gradient method to solve the reformulated optimization problem.
The effectiveness of the proposed algorithm is illustrated by several numerical examples in Section 5.

In the following, by $i\R$, we denote the imaginary axis of the complex plane, by $I_n$ the identity matrix of size $n \times n$, and by $\operatorname{null}(E)$ the nullspace of a matrix $E$.

\section{Dissipative Hamiltonian matrix pairs}\label{sec:DHpairs}

In this section, we construct  the setup to approach the nearest stable matrix pair problem.
As demonstrated in Example~\ref{ex:ex1}, the feasible set $\mathbb S$ is not open, not closed,
non-bounded, and highly nonconvex, thus it is very difficult to work directly with the set
$\mathbb S$. For this reason we  reformulate the nearest stable matrix pair problem into an equivalent optimization
problem with a simpler feasible set. 
For this, we extend the idea of a DH matrix from~\cite{GilS16} to dissipative Hamiltonian matrix pairs.
\begin{definition}
A matrix pair $(E,A)$, with $E,A \in \mathbb R^{n,n}$, is called a \emph{dissipative Hamiltonian (DH) matrix pair} if there exists an invertible matrix $Q\in \mathbb R^{n,n}$ such that
$Q^TE=E^TQ \succeq 0$, and $A$ can be expressed as $A=(J-R)Q$  with $J^T=-J$, $R^T=R  \succeq 0$.
\end{definition}
This definition of a DH matrix pair is a natural generalization of that of a DH matrix, because
for a standard DH matrix pair $(I_n,A)$ (when $E=I_n$),
$A$ is a DH matrix~\cite{GilS16}. Note however, that this  definition is slightly more restrictive than that of port-Hamiltonian descriptor systems in \cite{BeaMXZ17_ppt}, where it is not required that the matrix $Q$ is invertible.
In the following we call the matrix $R$ in a DH matrix pair $(E,(J-R)Q)$ the
\emph{dissipation matrix} as in the DH matrix case.

A DH-matrix pair is not necessarily regular as the following example shows.
\begin{example}\label{ex:ex2}{\rm
The pair $(E,(J-R)Q)$ with
\[
E=\begin{bmatrix}1 & 0&0\\0& 1&0\\0&0&0\end{bmatrix}, \;
J=\begin{bmatrix}0 & 2&0\\-2&0&0\\0&0&0\end{bmatrix}, \;
R=\begin{bmatrix} 1&0&0\\0&1&0\\0&0&0\end{bmatrix}, \;
Q=\begin{bmatrix}1 &0&0\\0 & 1&0\\0&0&1\end{bmatrix},
\]
is a DH matrix pair, but $(E,(J-R)Q)$ is singular, since
$\operatorname{det}(zE-(J-R)Q)\equiv 0$.
}
\end{example}
It is easy to see that
if the matrices $E$ and $A$ have a common nullspace,
then the pair $(E,A)$ is singular, but the converse is in general not  true, see e.~g. \cite{ByeHM98}.
However, for a singular DH matrix pairs, the property that $Q$ is invertible guarantees that the converse also holds.
\begin{lemma} \label{Lemma1}
Let $(E,(J-R)Q)$ be a DH-matrix pair. Then $(E,(J-R)Q)$ is singular if and only if
$\operatorname{null}(E)\cap\operatorname{null}((J-R)Q) \neq \emptyset$.
\end{lemma}
\proof The direction $(\Leftarrow)$  is immediate.  For the other direction, let $(E,(J-R)Q)$ be singular, i.~e. $\operatorname{det}(z E-(J-R)Q) \equiv 0$. Let $\lambda \in \mathbb C$ be such that $\real{(\lambda)} > 0$ and let  $x \in \mathbb C^n\setminus \{0\}$ be such that
\begin{equation}\label{eq:singular_DHpair_1}
(J-R)Qx=\lambda Ex.
\end{equation}
Since $Q$ is nonsingular, we have that $Qx\neq 0$, and we can multiply with $(Qx)^H$ from the left to obtain
\begin{equation}\label{eq:singular_DHpair_2}
x^HQ^TJQx-x^HQ^TRQx=\lambda x^HQ^TEx,
\end{equation}
where $x^H$ denotes the complex conjugate of  a vector $x$. This implies that $x^HQ^TEx=0$, because otherwise from~\eqref{eq:singular_DHpair_2} we would have
\[
\real{(\lambda)}=-\frac{x^HQ^TRQx}{x^*Q^TEx} \leq 0,
\]
since  $Q^TRQ \succeq 0$ (as $R \succeq 0$), and $Q^TE \succeq 0$. But this is a contradiction to the fact that
$\real{(\lambda)} > 0$. Therefore $x^HQ^TEx=0$ and also $Q^TEx=0$, and this implies that
$Ex =0$ as $Q$ is invertible. Inserting this in~\eqref{eq:singular_DHpair_1}, we get
$(J-R)Qx=0$, i.~e. $0\neq x\in \operatorname{null}(E)\cap\operatorname{null}((J-R)Q)$.
\eproof
Lemma~\ref{Lemma1} gives a necessary and sufficient condition for a DH-matrix pair to
be singular. However, if the dissipation matrix $R$ is positive definite, then regularity is assured, as shown in the following corollary.
\begin{corollary}\label{cor:DJ_regular}
Let $(E,(J-R)Q)$ be a DH matrix pair. If  $R$ is positive definite then the pair is regular.
\end{corollary}
\proof By Lemma~\ref{Lemma1}, a necessary condition for the pencil $\lambda E-(J-R)Q$
to be singular is that neither $E$ nor $(J-R)Q$ is invertible. Thus the result follows immediately
by the fact that if $R \succ 0$ in a DH-matrix pair $(E,(J-R)Q)$, then $(J-R)Q$ is invertible. Indeed,
suppose there exists $x\in \mathbb C\setminus \{0\}$ such that $(J-R)Qx=0$, then we have
$x^HQ^T(J-R)Qx=0$. This implies $x^HQ^TRQx=0$, since $J^T=-J$. 
Since  $Q$ is invertible and thus $Qx\neq 0$, this  is a contradiction to the assumption that $R$ is positive definite.
\eproof

In the following lemma, which is the matrix pair analogue of~\cite[Lemma 2]{GilS16},  we localize the finite eigenvalues of a DH-matrix pair.
\begin{lemma}\label{lem:DH_pair_prop}
Let $(E,(J-R)Q)$ be a regular DH-matrix pair and let $L(z):=zE-(J-R)Q$. Then the following statements hold.
\begin{enumerate}
\item All finite eigenvalues of the pencil $L(z)$ are in the closed left half of the complex plane.
\item The pencil $L(z)$ has a  finite eigenvalue $\lambda$ on the imaginary axis if and only if $RQx=0$ for some eigenvector
$x$ of the pencil $zE-JQ$ associated with  $\lambda$.
\end{enumerate}
\end{lemma}
\proof
Let $\lambda \in \mathbb C$ be an eigenvalue of $z E-(J-R)Q$ and let $x\in \mathbb C^n\setminus\{0\}$
be such that
\begin{equation}\label{eqn:DH_eig_1}
(J-R)Qx=\lambda Ex.
\end{equation}
Multiplying~\eqref{eqn:DH_eig_1} by $x^HQ^T$ from the left, we get
\begin{equation}\label{eqn:DH_eig_2}
x^HQ^T(J-R)Qx=\lambda\, x^HQ^TEx.
\end{equation}
Note that $x^HQ^TEx \neq 0$, because if $x^HQ^TEx=0$, then we have
$Q^TEx=0$ as $Q^TE \succeq 0$, and thus $Ex =0$ as $Q$ is invertible. Using this in~\eqref{eqn:DH_eig_1},
we have $(J-R)Qx=0$. This implies that $x \in \operatorname{null}(E)\cap \operatorname{null}((J-R)Q)$ which is,
by Lemma~\ref{Lemma1}, a contradiction to the regularity of the pair $(E,(J-R)Q)$.

Thus, by~\eqref{eqn:DH_eig_2} we have
\begin{equation}\label{eqn:DH_eig_3}
\real{(\lambda)}=-\frac{x^HQ^TRQx}{x^HQ^TEx} \leq 0,
\end{equation}
because $Q^TRQ \succeq 0$ as $R\succeq 0$, and $x^HQ^TEx >0$. This completes the proof of 1).

In the proof of 1), if  $\lambda \in i\mathbb R$, then from~\eqref{eqn:DH_eig_3} it follows that
$x^HQ^TRQx=0$. This implies that $RQx=0$, since $R\succeq 0$. Using this in~\eqref{eqn:DH_eig_1} implies that
$(\lambda E-JQ)x=0$.

Conversely, let $\lambda \in i\mathbb R$ and $x\in \mathbb C^n\setminus\{0\}$ be such that
$RQx=0$ and $(\lambda E-JQ)x=0$. Then this trivially implies that $\lambda$ is also an
eigenvalue of the pencil $\lambda E-(J-R)Q$ with eigenvector $x$. This completes the proof of 2).
\eproof
Making use of these preliminary results, we have the following stability characterization.
\begin{theorem}\label{thm:stable_semidef_R}
Every regular DH matrix pair $(E,(J-R)Q)$ of index at most one is stable.
\end{theorem}
\proof In view of Lemma~\ref{lem:DH_pair_prop}, to prove
the result it is sufficient to show that if $\lambda \in i\mathbb R$ is an
eigenvalue of the pencil $z E-(J-R)Q$, then $\lambda$ is semisimple.

Let us suppose that $\lambda \in i\mathbb R$ is a defective eigenvalue of the pencil $zE-(J-R)Q$
and the set $\{x_0,x_1,\ldots,x_{k-1}\}$ forms a Jordan chain of length $k$ associated with $\lambda$, see e.~g. \cite{GohLR82}, i.~e. $x_0 \neq 0$ and
\begin{eqnarray}
(\lambda E-(J-R)Q)x_0=0,\quad (\lambda E-(J-R)Q)x_1&=&Ex_0,\nonumber\\
(\lambda E-(J-R)Q)x_2&=&Ex_1,\nonumber\\
\vdots \label{thm:stable_DH_1}\\
 (\lambda E-(J-R)Q)x_{k-1}&=&Ex_{k-2}.\nonumber
\end{eqnarray}
Note that by Lemma~\ref{lem:DH_pair_prop}, we have that $(\lambda E-(J-R)Q)x_0=0$ implies that
\begin{equation}\label{thm:stable_DH_2}
(\lambda E-JQ)x_0=0\ \mbox{\rm and }\ RQ x_0=0.
\end{equation}
By~\eqref{thm:stable_DH_1}, $x_0$ and $x_1$ satisfy
\begin{equation}\label{thm:stable_DH_3}
(\lambda E-(J-R)Q)x_1=Ex_0.
\end{equation}
Multiplying~\eqref{thm:stable_DH_3} by $x_0^HQ^T$ from the left, we obtain
\[
x_0^H(\lambda Q^TE-(Q^TJQ-Q^TRQ))x_1=x_0^HQ^TEx_0.
\]
This implies that
\begin{eqnarray}\label{thm:stable_DH_4}
-x_1^H(\lambda Q^TE -Q^TJQ)x_0 +x_1^HQ^TRQ x_0 = x_0^HQ^TEx_0,
\end{eqnarray}
where the last equality follows by the fact that $Q^TE=E^TQ$, $J^T=-J$, and $R^T=R$.
Thus, by using~\eqref{thm:stable_DH_2} in~\eqref{thm:stable_DH_4}, we get
$x_0^HQ^TEx_0=0$. But this implies that $Q^TEx_0=0$ as $Q^TE \succeq 0$ and
$Ex_0 =0$ as $Q$ is invertible. Since $x_0$ is an eigenvector of the pencil $zE-(J-R)Q$  to
$\lambda$, $Ex_0=0$ implies that $(J-R)Qx_0=0$. This means that
$0\neq x_0  \in \operatorname{null}(E)\cap \operatorname{null}((J-R)Q)$, which contradicts the regularity
of the pair $(E,(J-R)Q)$. Therefore there does not exist a vector $x_1 \in \mathbb C^n$ satisfying~\eqref{thm:stable_DH_1}. Hence
$\lambda$ is semisimple.
\eproof
We note that the proofs of Lemma~\ref{lem:DH_pair_prop} and Theorem~\ref{thm:stable_semidef_R}
highly depend on the invertibility of $Q$. In fact, we have used the fact that
$(Q^TE,Q^T(J-R)Q)$ is regular, because $Q$ invertible implies that
$(E,(J-R)Q)$ is regular if and only if $(Q^TE,Q^T(J-R)Q)$ is regular.
If  $Q$ is singular and $Q^TE=E^TQ \succeq 0$, then the pair
$(Q^TE,Q^T(J-R)Q)$ is always singular, but the pair $(E,(J-R)Q)$ may be regular.

However, parts of Lemma~\ref{lem:DH_pair_prop} and Theorem~\ref{thm:stable_semidef_R}
also hold for singular $Q$ with an extra assumption on the eigenvectors of the
pencil $zE-(J-R)Q$, as stated in the next result. This result is
a generalization of~\cite[Lemma 2]{GilS16} for DH matrices with singular $Q$.

\begin{theorem}\label{eq:sing_Q}
Let $(E, (J-R)Q)$ be a regular matrix pair with $J^T=-J$, $R\succeq 0$, and
a singular $Q$ such that $E^TQ\succeq 0$. If no eigenvector $x$ of $L(z):=zE-(J-R)Q$
with respect to a nonzero eigenvalue satisfies $Q^TEx=0$, then
the following hold.
\begin{enumerate}
\item All finite eigenvalues of the pencil $L(z)$ are in the left half of the complex plane.
\item $L(z)$ has a finite eigenvalue $\lambda$ on the imaginary axis if and only if $RQx=0$ for some eigenvector
$x$ of the pencil $zE-JQ$ with respect to eigenvalue $\lambda$.
\item If $0\neq \lambda  \in i\mathbb R$ is an
eigenvalue of $z E-(J-R)Q$, then $\lambda$ is semisimple.
\end{enumerate}
\end{theorem}

Note that if $E=I_n$ in Theorem~\ref{eq:sing_Q}, then $Q$ is positive semidefinite and singular, and
the extra condition on the eigenvectors of pencil $zI_n - (J-R)Q$ (or, equivalently of matrix $(J-R)Q$)
is trivially satisfied. In this case Theorem~\ref{eq:sing_Q} coincides with the matrix case~\cite[Lemma 2]{GilS16}.

The following result from~\cite{MasKAS97} gives an equivalent condition for checking
asymptotical stability of such a pair $(E,A)$ by using Lyapunov's Theorem.
\begin{theorem}[\cite{MasKAS97}] \label{thm:impulsefree_reference}
Consider a pair $(E,A)$ with $E,A \in \mathbb R^{n,n}$. The pair is regular, of index at most one and asymptotically stable
if and only if there exists a nonsingular $V\in \mathbb R^{n,n}$ satisfying $V^TA+A^TV \prec 0$ and $E^TV=V^TE \succeq 0$.
\end{theorem}

\begin{theorem}\label{eq:main_result}
Let $(E,A)$ be a matrix pair, where $E,A \in \R^{n,n}$. Then the following statements are equivalent.
\begin{enumerate}
\item $(E,A)$ is a DH-matrix pair with positive definite dissipation matrix.
\item $(E,A)$ is regular, of index at most one, and asymptotically stable.
\end{enumerate}
\end{theorem}
\proof $1)\Rightarrow 2)$ Let $(E,A)$ be a DH-matrix pair with positive definite dissipation matrix, i.~e.,
$A=(J-R)Q$ for~some $R\succ 0$, $J^T=-J$, and nonsingular $Q$ with $Q^TE \succeq 0$.
Clearly, by Corollary~\ref{cor:DJ_regular} $(E,(J-R)Q)$ is regular. Furthermore,  $(E,(J-R)Q)$ has all its finite eigenvalues in the open left half plane. To see this,
 let $\lambda \in \mathbb C$ be a finite  eigenvalue of the pencil $zE-(J-R)Q$. Then by
Lemma~\ref{lem:DH_pair_prop} it follows that $\real{(\lambda)} \leq 0$, and $\real{(\lambda)} =0$
if and only if there exists $x\neq 0$ such that $(\lambda E-JQ)x=0$ and $0\neq Qx \in \operatorname{null}(R)$. But
$\operatorname{null}(R) =\{0\}$ as $R \succ 0$. 

To show that $(E,(J-R)Q)$ is of index at most one, we
set $r:=\text{rank}(E)$ and assume that $U \in \R^{n,n-r}$ is an orthogonal matrix whose column spans $\text{null}(E)$.
Then, see \cite{KauNC89}, $(E,(J-R)Q)$ is of index at most one if and only if $\text{rank}(\mat{cc}E&(J-R)QU\rix)=n$.
Suppose that $x \in \C^n\in \setminus \{0\}$ is such that $x^H\mat{cc}E&(J-R)QU\rix =0$. Then we have the two conditions
\begin{equation}\label{eq:equivalent_proof_1}
x^HE=0,\ x^H(J-R)QU=0.
\end{equation}
Since $Q$ is invertible, we have $x^HEQ^{-1}=0$ and hence
$(EQ^{-1})x=0$ because $EQ^{-1} \succeq 0$ as $E^TQ \succeq 0$. This shows that $Q^{-1}x \in \text{null}(E)$,
and thus there exists $y \in \C^{n-r}$ such that $Q^{-1}x=Uy$, or, equivalently $x=QUy$.
Using this in~\eqref{eq:equivalent_proof_1}, we obtain that $x^H(J-R)x=0$. This implies that
$x^HJx=0$ and $x^HRx=0$ as $J$ is skew-symmetric and $R$ is symmetric. But this is a contradiction to the assumption that $R \succ 0$. This completes the proof of $1)\Rightarrow 2)$.

$2)\Rightarrow 1)$ Consider a  pair $(E,A)$, with $E,A\in \mathbb R^{n,n}$, that is regular, asymptotically stable, and of index at most one. Then by
Theorem~\ref{thm:impulsefree_reference}, there exist an nonsingular $V\in \mathbb R^{n,n}$ such that
$V^TA+A^TV \prec 0$ and $\quad E^TV=V^TE \succeq 0$.
Setting
\[
Q=V,\quad J=\frac{(AV^{-1})-(AV^{-1})^T}{2},\quad \text{and}\quad R=-\frac{(AV^{-1})+(AV^{-1})^T}{2},
\]
we have $J^T=-J$, $E^TQ=Q^TE \succeq 0$, and $R \succ 0$, as $V$ is invertible. Applying the Lyapunov inequality
\[
V^TRV=-\frac{V^T((AV^{-1})+(AV^{-1})^T)V}{2}=-\frac{V^TA+A^TV}{2} \succ 0,
\]
the assertion follows.
\eproof
We conclude this section with the observation that the set of DH matrix pairs is invariant under orthogonal transformations of the
matrix pair.
\begin{lemma}\label{lem:ortho_invariant}
Let $(E,A)$ be a DH matrix pair, and let $U$ and $V$ be orthogonal matrices such that
$\widetilde E=UEV$ and $\widetilde A=UAV$. Then $(\widetilde E,\widetilde A)$ is a DH pair.
\end{lemma}
\proof
Since $(E,A)$ is DH, we have $A=(J-R)Q$ for some $J^T=-J$, $R\succeq 0$, and invertible $Q$
such that $Q^TE \succeq 0$. Using the orthogonality of $U$ and $V$, we have
\begin{eqnarray*}
(\widetilde E,\widetilde A)&=& (UEV,UAV)=(UEV,U(J-R)QV)=(UEV,U(J-R)U^TUQV)\\
&=& (UEV,(UJU^T-URU^T)UQV) = (\widetilde E-(\widetilde J-\widetilde R)\widetilde Q),
\end{eqnarray*}
where $\widetilde J:=UJU^T$, $\widetilde R:=URU^T$, and $\widetilde Q:=UQV$. Since
$U$ and $V$ are orthogonal, we have  that $\widetilde R\succeq 0$ as $R \succeq 0$, and $\widetilde Q$ is invertible with
$\widetilde Q^T \widetilde E \succeq 0$ as $Q$ is invertible with $Q^TE \succeq 0$. This shows that
$(\widetilde E,\widetilde A)$ is again DH.
\eproof

\section{Reformulation of the nearest stable matrix pair problem}

In this section, we exploit the results obtained in the previous section and derive a
reformulation of the nearest stable matrix pair problem.
By Theorem~\ref{eq:main_result},
the set $\mathbb S$ of all asymptotically stable matrix pairs can be expressed as the set of all DH matrix pairs with positive definite dissipation, i.~e.,
\begin{equation*}
\mathbb S=\left\{(E,(J-R)Q)\in \R^{n,n}\times \R^{n,n}~:~J^T=-J,\ R\succ 0,\ Q~\text{invertible s.t.}~Q^TE \succeq 0\right\}=:\mathbb S_{DH}^{\succ 0}.
\end{equation*}
This characterization changes the feasible set and also the objective function in the nearest stable matrix pair problem as
\begin{equation*}
\inf_{(M,X) \in \mathbb S} \{{\|E-M\|}_F^2+{\|A-X\|}_F^2\}=
\inf_{(M,(J-R)Q) \in \mathbb S_{DH}^{\succ 0}} \{{\|E-M\|}_F^2+{\|A-(J-R)Q\|}_F^2\}.
\end{equation*}
We have demonstrated in Example~\ref{ex:ex1} that the set $\mathbb S$ of all stable matrix pairs is neither open nor
closed and clearly the alternative characterization of  $\mathbb S$
in terms of $\mathbb S_{DH}^{\succ 0}$ does not change this, since $\mathbb S_{DH}^{\succ 0}$ is not closed due to the
constraint that $R \succ 0$ and $Q$ is invertible, and not open due to the constraint  $Q^TE \succeq 0$.
However, we can instead look at the set $\mathbb S_{DH}^{\succeq 0}$ containing all pairs of the form $(E,(J-R)Q)$ with
$J^T=-J$, $R \succeq 0$ ($R$ can be singular), and $Q$ ($Q$ can be singular) such that $E^TQ \succeq 0$, i.~e.
\[
\mathbb S_{DH}^{\succeq 0}:=\left\{(E,(J-R)Q)\in \R^{n,n}\times \R^{n,n}~:~J^T=-J,\ R\succeq 0,\ Q^TE \succeq 0\right\}.
\]
Then
$\mathbb S_{DH}^{\succeq 0}$ is the closure $\overline{\mathbb S_{DH}^{\succ 0}}$ of $\mathbb S_{DH}^{\succ 0}$.
Following the arguments similar to that of Lemma~\ref{lem:ortho_invariant} we also have that $\mathbb S_{DH}^{\succeq 0}$ is
invariant under orthogonal transformations and, furthermore, we have that
\begin{equation}\label{eq:reformulation_1}
\inf_{(M,(J-R)Q) \in \mathbb S_{DH}^{\succ 0}} \{{\|E-M\|}_F^2+{\|A-(J-R)Q\|}_F^2\}=
\inf_{(M,(J-R)Q) \in \mathbb S_{DH}^{\succeq 0}} \{{\|E-M\|}_F^2+{\|A-(J-R)Q\|}_F^2\}.
\end{equation}
Note, however, that the set $\mathbb S_{DH}^{\succeq 0}$ is  not bounded, and hence the infimum in the right hand side of~\eqref{eq:reformulation_1} may not be attained. 

Our observations lead to the following  reformulation of the nearest stable matrix pair problem.
\begin{theorem}\label{thm:reformulation_alg}
Let $(E,A)\in \R^{n,n}\times \R^{n,n}$. Then
\begin{equation}\label{thm_eq:reform_1}
\inf_{(M,X) \in \mathbb S} \{{\|E-M\|}_F^2+{\|A-X\|}_F^2\}=
\inf_{(M,(J-R)Q) \in \mathbb S_{DH}^{\succeq 0}} \{{\|E-M\|}_F^2+{\|A-(J-R)Q\|}_F^2\}.
\end{equation}
\end{theorem}
For the standard system (when $E$ is the identity matrix in~\eqref{eq:lindae}) stability solely depends on
the eigenvalues of the matrix $A$. Thus making $A$ stable without perturbing the identity matrix gives an
upper bound for the distance of $(I_n,A)$ to the nearest stable matrix pair.
This also follows from~\eqref{thm_eq:reform_1} because
\begin{equation}\label{thm_eq:reform_2}
\inf_{(M,(J-R)Q) \in \mathbb S_{DH}^{\succeq 0}} \{{\|I_n-M\|}_F^2+{\|A-(J-R)Q\|}_F^2\}
\leq \inf_{(I_n,(J-R)Q) \in \mathbb S_{DH}^{\succeq 0}} \{{\|A-(J-R)Q\|}_F^2\}.
\end{equation}
We note that the infimum on the right hand side of~\eqref{thm_eq:reform_2}
is the distance of $A$ from the set of stable matrices~\cite{GilS16}. We will demonstrate in our numerical experiments that (as  expected) the inequality in~\eqref{thm_eq:reform_2} is strict.
However, it is interesting to note that if $(J^*,R^*,Q^*)$ is a stationary point of the right hand side of~\eqref{thm_eq:reform_2} then $(I_n,J^*,R^*,Q^*)$ is a stationary point of the left hand side of~\eqref{thm_eq:reform_2}; see the discussion in Section~\ref{algosstable}.

A DAE with coefficient pair $(E,A)$ and nonsingular $E$ can be equivalently reformulated as a standard system
$(I_n,AE^{-1})$, and then stability can be determined by the eigenvalues of
 $AE^{-1}$. Thus again making $AE^{-1}$ stable, or equivalently, $(I_n,AE^{-1})$
stable without perturbing the identity matrix, gives an upper bound
for the distance of $(E,A)$ to the nearest stable matrix pair. Indeed,
w.l.o.g.\@ if we scale $(E,A)$ so that  ${\|E\|}_F=1$,
then
\begin{eqnarray}
\nonumber&&\inf_{(M,(J-R)Q) \in \mathbb S_{DH}^{\succeq 0}} \{{\|E-M\|}_F^2+{\|A-(J-R)Q\|}_F^2\}\\
&&\qquad\leq \inf_{(E,(J-R)Q) \in \mathbb S_{DH}^{\succeq 0}} \{{\|A-(J-R)Q\|}_F^2\} \nonumber\\
&& \qquad= \inf_{(I_n,(J-R)P) \in \mathbb S_{DH}^{\succeq 0}} \{{\|A-(J-R)PE\|}_F^2\} \label{thm_eq:reform_3}\\
&&\qquad\leq \inf_{(I_n,(J-R)P) \in \mathbb S_{DH}^{\succeq 0}} \{{\|AE^{-1}-(J-R)P\|}_F^2\} \label{thm_eq:reform_4},
\end{eqnarray}
where~\eqref{thm_eq:reform_3} follows by first using the fact that
$Q^TE=E^TQ \succeq 0$ implies that $E^{-T}Q^T=QE^{-1} \succeq 0$,
whereas in~\eqref{thm_eq:reform_4} we use the sub-multiplicativity of the
Frobenius norm and ${\|E\|}_F=1$. Thus, for nonsingular $E$,
an upper bound can be computed by the reduction to a standard system. However,
the reduction process may be numerically unstable. Hence, it is
advisable to work directly with the pair $(E,A)$ even for regular
matrix pairs with invertible $E$.
\begin{remark}\label{rem:sensitivity}{\rm
The eigenvalues of the regular DH pair $(E,(J-R)Q)$ can be
highly sensitive to small perturbations. This happens  if the finite or infinite eigenvalues are multiple
and defective or if they are close to being multiple and the eigenvectors have a small angle.
In fact, let $\lambda$ be a simple eigenvalue
of $(E,(J-R)Q)$ and let $x$ be the corresponding eigenvector normalized to unit norm.
By~\eqref{eqn:DH_eig_3} we have
\[
\real{(\lambda)} = -\frac{x^*Q^TRQx}{x^*Q^TEx},
\]
which implies that $(-\real{(\lambda)},x)$ is an eigenvalue/eigenvector pair of the symmetric pencil
$L(z):=zQ^TE-Q^TRQ$, since $Q^TE \succeq 0$ and $Q^TRQ \succeq 0$.
Then, see e.~g.~\cite{Tis00}, the normwise condition number of $\lambda$  of $L(z)$
is given by
\[
\kappa(-\real{(\lambda)},L) = \frac{1+|\real{(\lambda)}|^2}{x^*Q^TEx},
\]
which implies that
\[
\real{(\lambda)}=-\frac{x^*Q^TRQx}{x^*Q^TEx}=-\frac{x^*Q^TRQx}{1+|\real{(\lambda)}|^2}\,\kappa{(-\real{(\lambda)},L)}.
\]
Therefore, if $\kappa(-\real{(\lambda)},L)$ is large, then a small perturbation
can significantly perturb the eigenvalues of $(E,(J-R)Q)$.
}
\end{remark}

\section{Algorithmic solution to the nearest stable matrix pair problem} \label{algosstable}
In this section, we propose an algorithm to solve the nearest stable matrix pair problem
using the reformulation of Theorem~\ref{thm:reformulation_alg}, i.~e., to solve
\begin{equation} \label{originaloptprob}
\min_{J = -J^T,\ R \succeq 0,\ Q,\ M,\ Q^T M \succeq 0}  \{{\|E - M\|}_F^2 + {\| A - (J-R)Q \|}_F^2\}.
\end{equation}
One of the first algorithms that comes to mind to solve~\eqref{originaloptprob} is a block-coordinate descent (BCD) method (see for example \cite{wright2015coordinate}) with three blocks of variables $(J,R)$, $Q$ and $M$.
In this case
the subproblems for each of the three blocks are convex when the others are fixed (these are least-squares problems with
linear or positive semidefinite constraints).
However, we have observed that this method is not very efficient in practice, since convergence is slow when one gets close to a stationary point. Moreover,  the BCD method can get easily stuck in saddle points. In particular, consider a matrix pair of the form $(I_n,A)$ and let $(J,R,Q)$ be a fixpoint of the BCD method (with two blocks of variables $(J,R)$ and $Q$) for the nearest stable matrix problem
\[
\min_{J = -J^T, R \succeq 0, Q \succeq 0} \{{\| A - (J-R)Q \|}_F^2\} .
\]
Then $(I_n,J,R,Q)$ is a fixpoint of the BCD method (with three blocks of variables $(J,R)$, $Q$ and $M$)
for~\eqref{originaloptprob}. By construction, then the update of $(J,R)$ and $Q$ cannot be improved
for $M=I_n$ fixed, while $M=I_n$ is optimal, since it is feasible and minimizes
${\|M-I_n\|}_F$. This behavior of the BCD method is illustrated in the following example.
\begin{example}\label{ex1} {\rm
Consider a matrix pair $(E,A)$, where
\[
A = \begin{bmatrix}1 & 1 & 0 \\ -1& 1&1\\0&-1&1\end{bmatrix}\quad \text{and}\quad E = I_n,
\]
 with eigenvalues $1$, $1\pm \sqrt{2}i$.
The representation $(J,R,Q,I_n)$ with
\[
J = \begin{bmatrix}0 & 1 & 0 \\ -1& 0&1\\0&-1&0\end{bmatrix},~R = 0,~ Q = I_n,~\text{and}~ M = I_n,
\]
is a fixpoint of the BCD method for~\eqref{originaloptprob} with error ${\|A-(J-R)Q\|}_F^2=3$ and we strongly believe that $(J-R)Q$ is the nearest stable matrix to $A$.

Using our fast projected gradient method that will be introduced below, initialized with this solution, we obtain a nearest stable matrix pair  with much lower distance
\[
{\|A-(J-R)Q\|}_F^2 + {\|M-I_n\|}_F^2=1.536.
\]
}
\end{example}

One possible reason why the formulation~\eqref{originaloptprob} does not seem to lead to good solutions is that the constraint $Q^T M \succeq 0$ couples the variables $Q$ and $M$ rather strongly. This motivated us to introduce a further reformulation of~\eqref{originaloptprob}, where the feasible set is free of any constraint that involves coupling of some of its variables
and onto which it is easy to project.

\subsection{ Reformulation of~\eqref{originaloptprob} }

The formulation~\eqref{originaloptprob} is not convenient for standard optimization schemes, such as a projected gradient method, because the feasible set is highly non-convex with the constraints $Q^T M \succeq 0$, and this constraint seems to prevent standard optimization schemes to converge to good solutions, as demonstrated in Example~\ref{ex1}. To overcome this dilemma
we introduce a new variable  $H = Q^TM = M^T Q \succeq 0$ in~\eqref{originaloptprob}, which corresponds to the Hamiltonian for port-Hamiltonian descriptor systems, see \cite{BeaMXZ17_ppt}.  This results in another reformulation of~\eqref{originaloptprob} with a modified  feasible set and objective function
\begin{equation} \label{reformoptprob}
\min_{J = -J^T, R \succeq 0, Q, H \succeq 0} {\|A - (J-R) Q\|}_F^2 + {\|E - Q^{-T}H\|}_F^2 .
\end{equation}
The feasible set of~\eqref{reformoptprob} is rather simple, with no coupling of the variables,
and it is relatively easy to project onto it.

\begin{remark}{\rm
If one want to obtain a DH pair with $Q$ invertible and $R \succ 0$ so that $(M,(J-R)Q)$
is regular, asymptotically stable and is of index at most one, see Theorem~\ref{eq:main_result}, then
the constraints $R \succeq 0$ and $H \succeq 0$ can be replaced with
$R \succeq \delta I_n$ and $H \succeq  \delta I_n$ for some small parameter $\delta > 0$ that prevents $R$ and $H =Q^T M$ from being
singular. As a consequence then $Q$ will be invertible, since otherwise $Q^{-T}H$ is unbounded.
Numerically, it does not make the problem much different, the projection is still rather straightforward.
Furthermore, in view of Remark~\ref{rem:sensitivity}, the constraints $R \succeq \delta I_n$ and $H \succeq  \delta I_n$ will assure that every finite eigenvalue $\lambda$ of the DH pair $(M,(J-R)Q)$ satisfies
\[
\real{(\lambda)}=-\frac{x^*Q^TRQx}{x^*Q^TMx}=-\frac{x^*Q^TRQx}{1+|\real{(\lambda)}|^2}\,\kappa{(-\real{(\lambda)},L)}
\leq -\frac{\delta}{1+|\real{(\lambda)}|^2}\,\kappa{(-\real{(\lambda)},L)},
\]
where $\kappa{(-\real{(\lambda)},L)}$ is the condition number of $-\real{(\lambda)}$ being a simple eigenvalue
 of the symmetric pencil $L(z)=zQ^TM-Q^TRQ$.
}
 \end{remark}

\subsection{Fast projected gradient algorithm} \label{subsecalgo}

Following the spirit of~\cite{GilS16}, we use a fast projected gradient algorithm to solve~\eqref{reformoptprob} despite the fact that it is a non-convex problem. The only non-trivial part in implementing this method is to compute the gradient of the objective function with respect to $Q$. We show in Appendix~\ref{appgradQ} that
\[
\frac{1}{2} \nabla_Q f(J,R,Q,H)
= (J-R)^T ( (J-R)Q - A )
+ Q^{-T}H (E^T - HQ^{-1}) Q^{-T},
\]
where $f(J,R,Q,H) := {\|A - (J-R) Q\|}_F^2 + {\|E - Q^{-T}H\|}_F^2$.
A pseudocode for our fast projected gradient methods is presented in Algorithm~\ref{fastgrad}.

 \algsetup{indent=2em}
\begin{algorithm}[ht!]
\caption{Fast Gradient Method (FGM)~\cite[p.90]{Nes04}} \label{fastgrad}
\begin{algorithmic}[1]
\REQUIRE The (non-convex) function $f(x)$, an initial guess $x \in \mathcal{X}$,
a number of iterations $K$ 

\ENSURE An approximate solution $x \approx \argmin_{z \in \mathcal{X}} f(z)$.  \medskip

\STATE $\alpha_1 \in (0,1)$; $y = x$ ; initial step length $\gamma > 0$.

\FOR{$k = 1 :$ $K$}

\STATE $\hat x = x$. \hspace{2.55cm} \emph{\% Keep the previous iterate in memory.}

\STATE $x = \mathcal{P}_{\mathcal{X}}  \left(  y - \gamma \nabla f(y) \right)$. \emph{\% $\mathcal{P}_{\mathcal{X}}$ is the projection on $\mathcal{X}$}

\WHILE { $f(x)  > f(\hat x)$ and $ \gamma > \underline{\gamma}$ }

\STATE Reduce $\gamma$.

\STATE $x = \mathcal{P}_{\mathcal{X}}  \left(  y - \gamma \nabla f(y) \right)$.

\ENDWHILE

\IF { $\gamma < \underline{\gamma}$ }

\STATE Restart fast gradient: $y = x$; $\alpha_{k} = \alpha_{0}$.

\ELSE

\STATE $y = x + \beta_k \left(x - \hat x\right)$, \quad where $\beta_k =  \frac{\alpha_{k} (1-\alpha_{k})}{\alpha_{k}^2 + \alpha_{k+1}}$
with
$\alpha_{k+1} \geq 0$ s.t.
$\alpha_{k+1}^2 = (1-\alpha_{k+1}) \alpha_{k}^2$.

\ENDIF

\STATE Increase $\gamma$.

\ENDFOR

\end{algorithmic}
\end{algorithm}

\paragraph{Initialization} To  initialize the fast projected gradient approach we use a similar strategy as in the matrix case, see~\cite{GilS16}.  For $Q = I_n$, the optimal solutions for $J$, $R$ and $H$ in~\eqref{reformoptprob} can be computed easily, and are given by
\[
J = \frac{A-A^T}{2},
R = \mathcal{P}_{\succeq }\left(\frac{-A-A^T}{2}\right),
H = \mathcal{P}_{\succeq }\left(E^T\right),
\]
where $\mathcal{P}_{\succeq }$ is the projection onto the set of positive semidefinite matrices.

Because of the non-convex nature of the problem, the solution obtained by our algorithm is highly sensitive to the initial
point. We observed that the above initial point provides in general good solutions. Future work will require the development or more sophisticated initialization schemes and possibly globalization approaches to obtain better solutions.

\paragraph{Line-search}

For $Q$ fixed, the Lipschitz constant of the gradient of $f$ with respect to
$H$ (resp.\@ $(J,R)$) is given by $L=\lambda_{\min}(QQ^T)^{-1}$ (resp.\@ $L=\lambda_{\max}(QQ^T)$),
where $\lambda_{\min}(X)$ (resp.\@ $\lambda_{\max}(X)$) denotes the smallest (resp.\@ largest) eigenvalue of matrix $X$.
Hence, performing a gradient step with step length $\gamma = 1$ in the variable $H$ (resp.\@ $(J,R)$) will guarantee the decrease of the objective function. Since we initialize with $Q = I_n$, we choose an initial step length of $\gamma = 1$.
If it does not lead to a decrease of the objective function, then we decrease $\gamma$ by a constant factor (we used 2/3).
This is a standard backtracking line-search strategy. If no decrease of the objective function is possible, i.~e. if the step length is smaller than some threshold $\underline{\gamma}$, then we restart the fast gradient scheme with a standard gradient step, which guarantees convergence.
At the next step, we use a slightly larger step length than the previous step to avoid the step length to go to zero, e.~g.\@ by multiplying $\gamma$ by a constant factor; we used 2.

\section{ Numerical Examples }

Our optimization code is available from \url{https://sites.google.com/site/nicolasgillis/} and the numerical examples presented below can be directly run from this online code.
The algorithm runs in $O(n^3)$ operations, including projections on the set of positive semidefinite matrices, inversion of the matrix $Q$ and all necessary matrix-matrix products. Our algorithm can be applied on a standard desktop computer with $n$ up to a thousand.
As far as we know, no other algorithm has been proposed for the computation of the nearest stable matrix pair. Therefore, we cannot compare to an existing method, so for illustration of our results we will only compare our fast gradient method with the projected gradient method (GM), which  is simply FGM where restart is used at every step. For all experiments, we limit the CPU time to 10 seconds.

\subsection{ Case 1: $(I_n,A)$ } \label{case1}

Let us take as $A$ the Grcar matrix from~\cite{GilS16} and $E = I_n$. The Grcar matrix of order $k$ is a
banded Toeplitz matrix with its subdiagonal set to $-1$, and both its main and $k$ superdiagonals set to 1.

\begin{figure*}[h!]
\centering
\begin{tabular}{cc}
\includegraphics[width=0.5\textwidth]{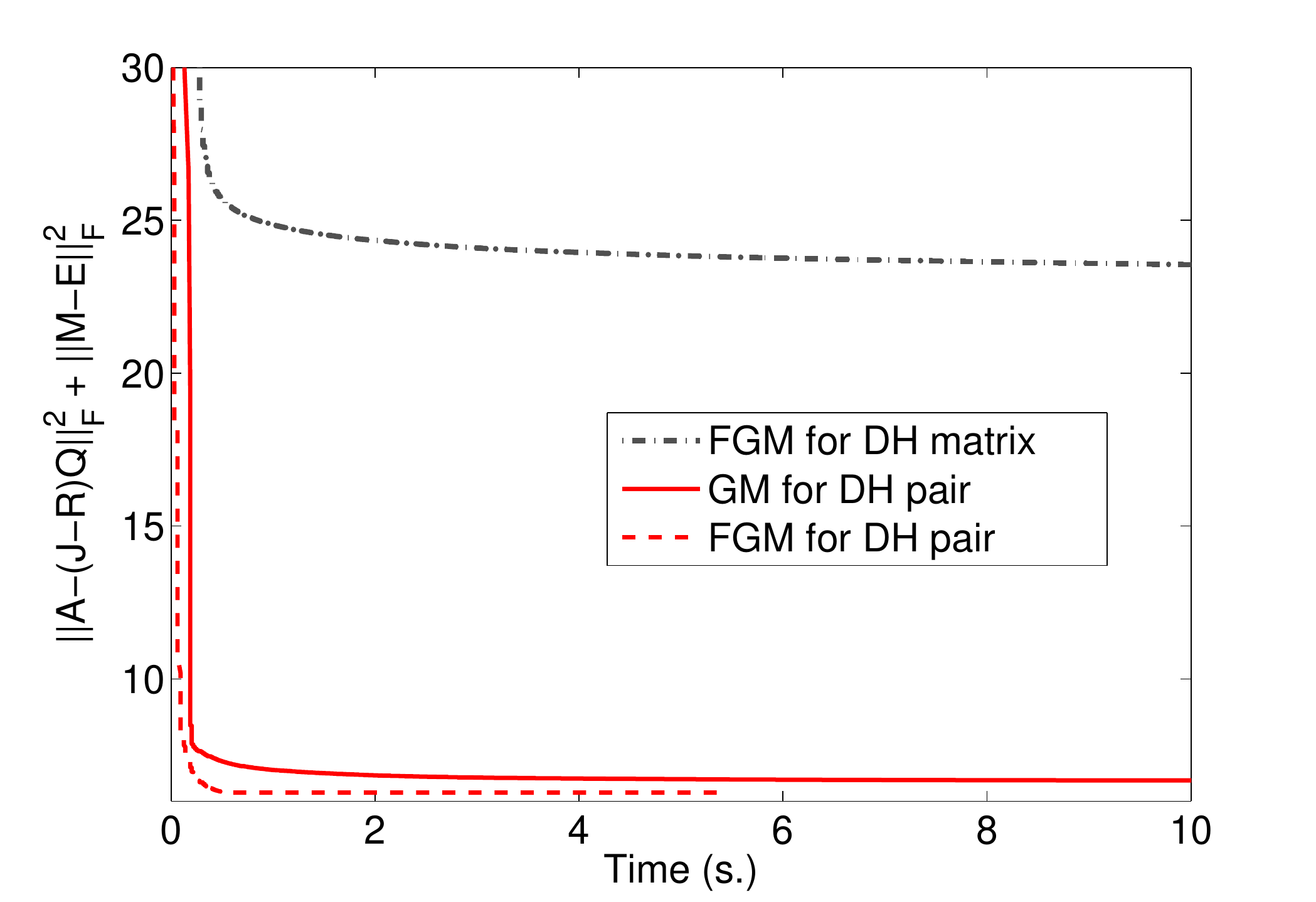} &  \includegraphics[width=0.5\textwidth]{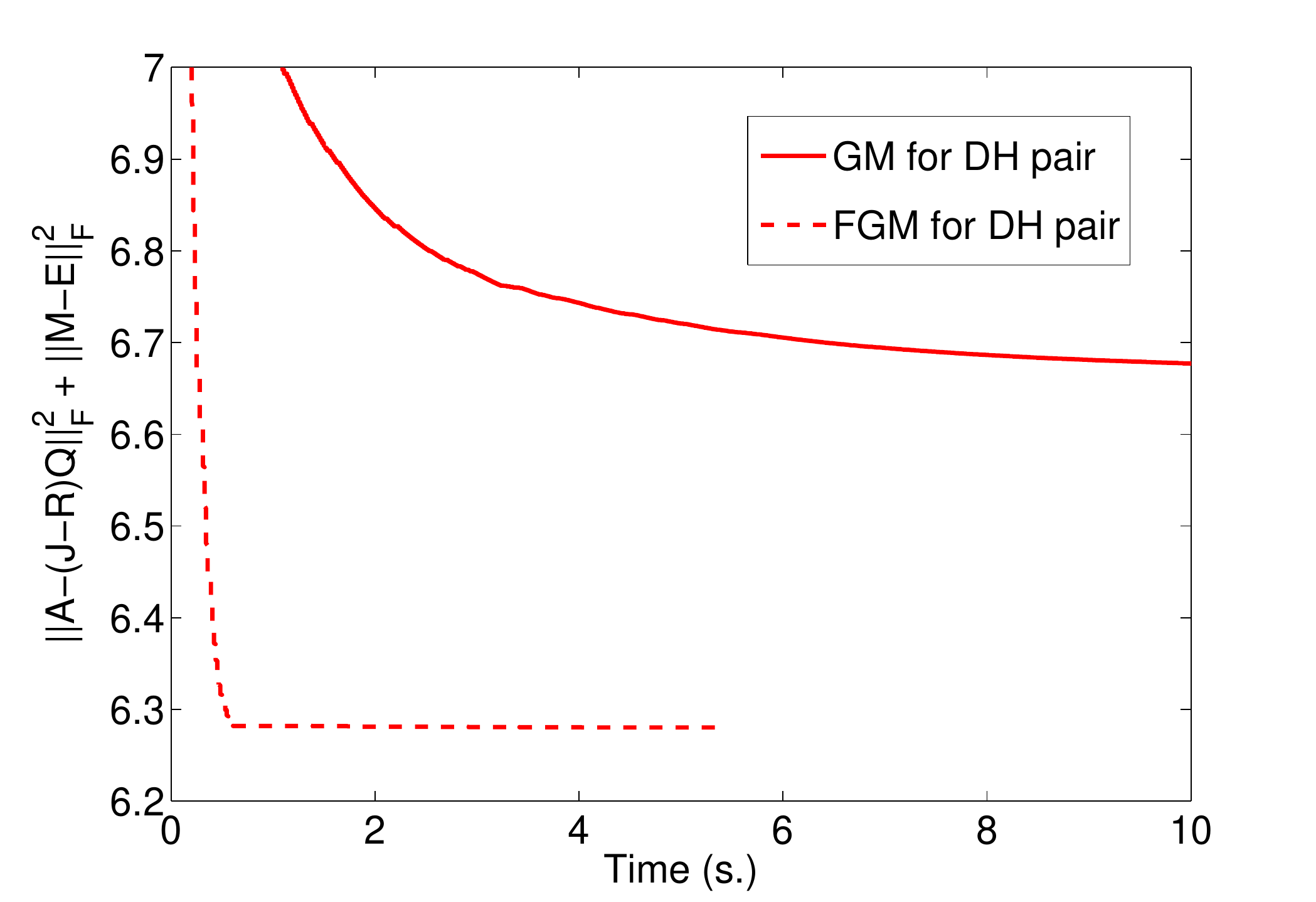}
\end{tabular}
\caption{Evolution of the objective function for FGM and GM for the matrix pair $(I_n,A)$ where $A$ is a $20\times 20$ Grcar matrix of order $3$. The right figure is a zoom of the left figure and shows the fast convergence of FGM compared to GM.
\label{grcar203}}
\end{figure*}

We use $n=20$ and $k=3$.
The nearest stable matrix $\tilde{A}$ to $A$ found in~\cite{GilS16} satisfies ${\|A-\tilde{A}\|}_F^2 = 23.51$.
Figure~\ref{grcar203} shows the evolution of the objective function of~\eqref{originaloptprob}.
Our FGM achieves a objective function value ${\|A - (J-R) Q\|}_F^2 + {\|E - Q^{-T}H\|}_F^2 = 6.28$.
Note that FGM converges much faster than GM.

\subsection{ Case 2: Random $A$ and rank-deficient $E$ }

In our second example, we generate each entry of $A$ using the normal distribution of mean $0$ and standard deviation $1$ (\texttt{randn(n)} in Matlab). We generate $E$ so as it has rank $r$, by taking $E$ as the best rank-$r$ approximation of a matrix generated randomly in the same way as $A$. We use $n = 20$ and $r = 3$.

For these types of matrices, usually the algorithm converges to a rank-deficient solution, i.~e. $R,H,Q$ are not of full rank,
hence it is useful to use the lower bound for the eigenvalues to obtain better numerical stability; we use $\delta =10^{-6}$.

\begin{figure}[h!]
\centering
\includegraphics[width=0.5\textwidth]{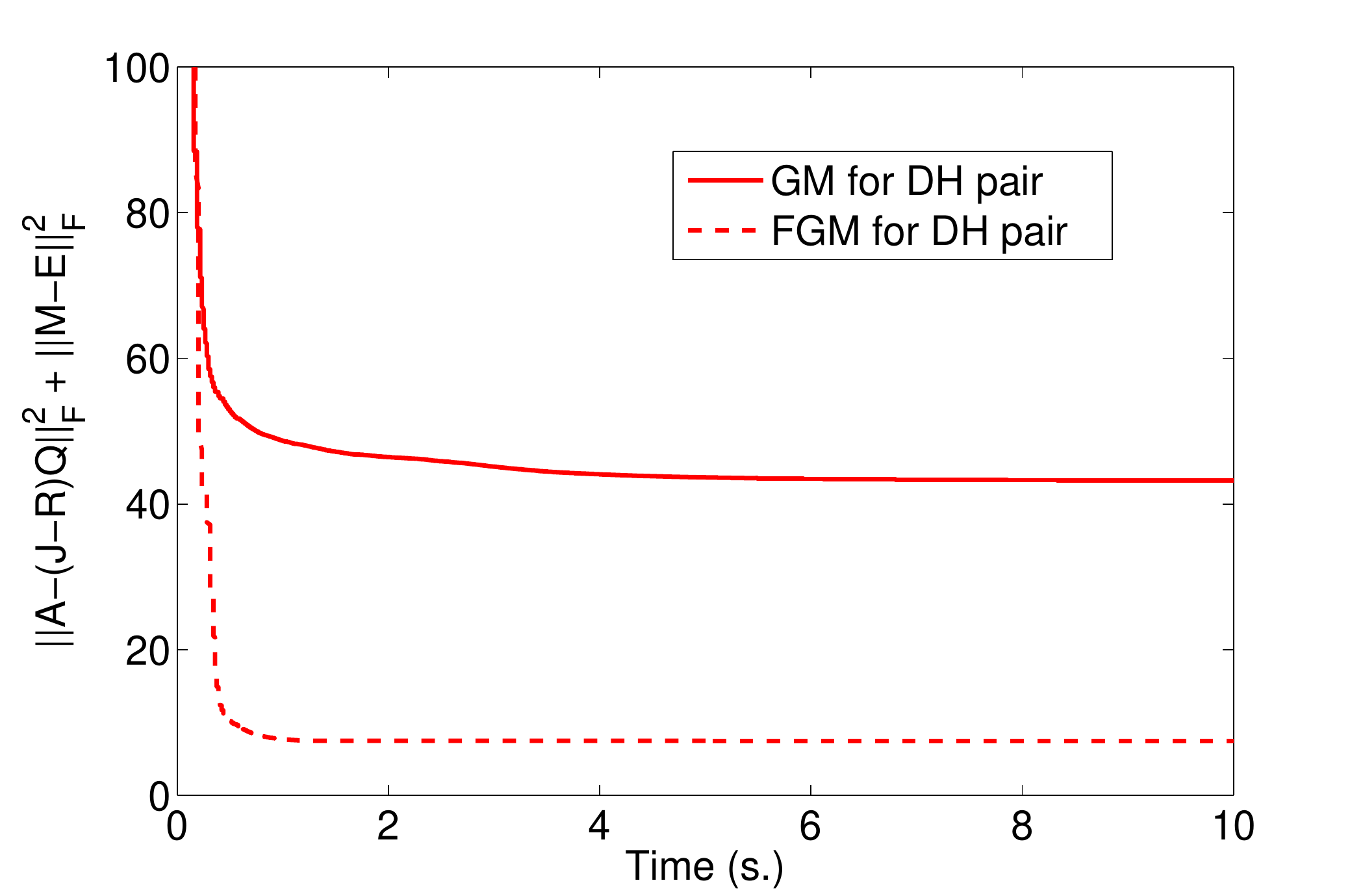}
\caption{Evolution of the objective function for the different algorithms for matrix pair $(E,A)$ where $A$
is Gaussian and $E$ is the best rank-$r$ approximation of a Gaussian matrix, $n=20$ and $r=3$.
\label{rand203}}
\end{figure}

Figure~\ref{rand203} shows the evolution of the objective function of~\eqref{originaloptprob} for a particular example.
On all the examples that we have run, FGM always converged much faster and generated a significantly better solution than GM, similarly as shown on Figure~\ref{rand203}.

\subsection{ Case 3: Mass-spring-damper model}

Let us consider a lumped parameter, mass-spring-damper dynamical system,
see, e.~g.,~\cite{Ves11} with $n$ point masses and $n$ spring-damper pairs. The equations
of motion can be written in the form of second order system of differential equations as
\begin{equation}\label{eq:MSD}
M\ddot{u}+D\dot{u}+Ku=f,
\end{equation}
where $M \succ 0$ is the mass matrix, $D \succeq 0$ is the damping matrix and
$K \succ 0$ is the stiffness matrix, $f,u,\dot{u},\ddot{u}$ are the forces, displacement,
velocity and acceleration vectors, respectively.

An appropriate first order formulation of~\eqref{eq:MSD} is associated with the DH pencil
$\lambda E-A$, where
\begin{equation} \label{msdmat}
E =  \mat{cc} M & 0\\ 0  &  I_{n} \rix,~
A = (J-R) Q,~
J =   \mat{cc} 0 & -I_{n} \\ I_{n} & 0 \rix ,~
R =  \mat{cc} D & 0 \\ 0 & 0 \rix ,~
Q =  \mat{cc} I_{n} & 0 \\ 0 & K \rix.
\end{equation}
Since $E$ is invertible, the pair $(E,A)$ is regular and of index zero, and thus from Theorem~\ref{thm:stable_semidef_R}
is also stable. In order to make the pair unstable, we perturb the dissipation matrix to become indefinite via
\[
R =  \mat{cc} D & 0 \\ 0 & -\epsilon I_n \rix.
\]
We use $n = 10$ and $\epsilon = 0.1$ which moves some eigenvalues to the open right half of the complex plane (see Figure~\ref{MSD101toneig} below).

\begin{figure*}[h!]
\centering
\includegraphics[width=0.6\textwidth]{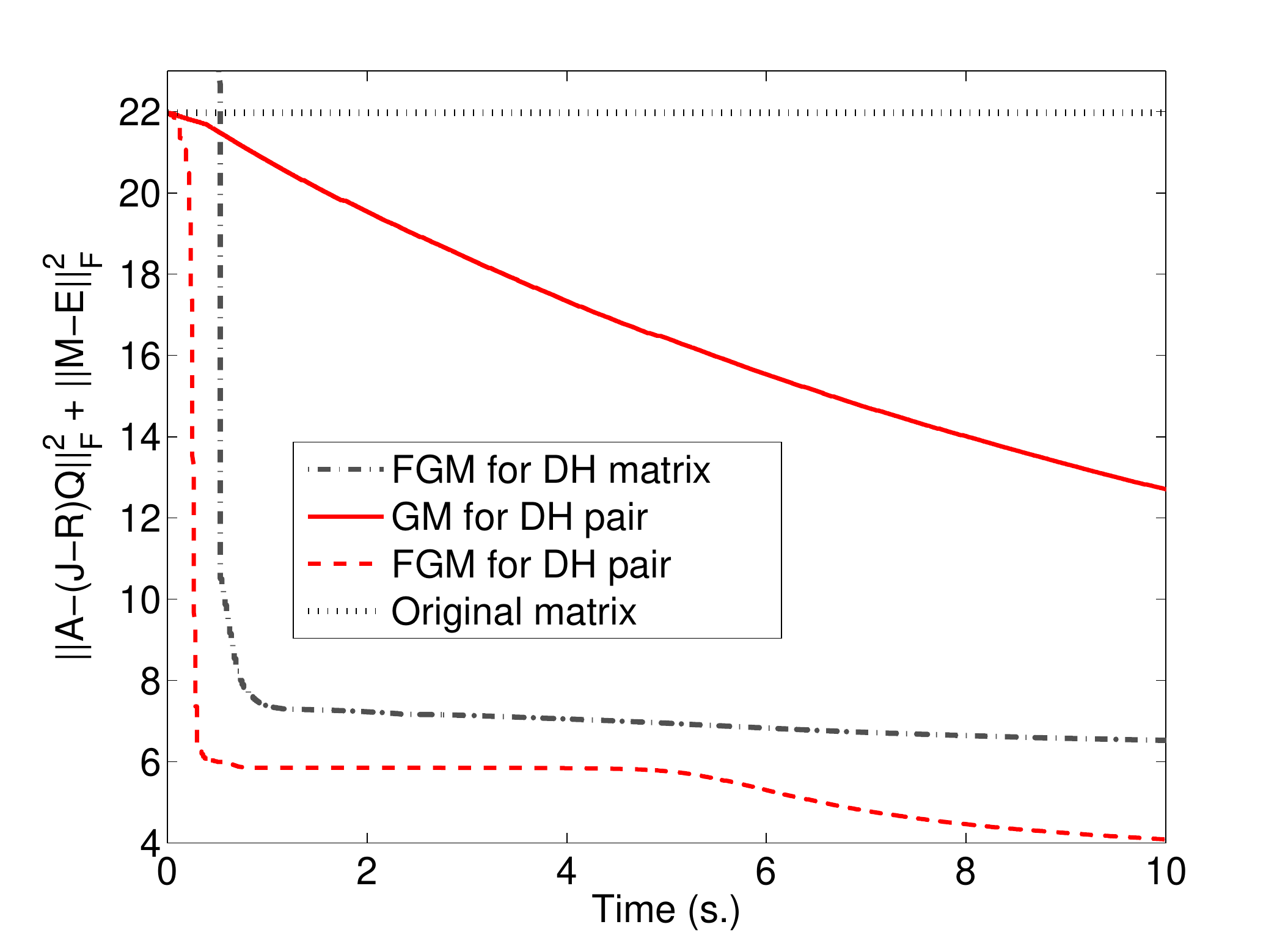}
\caption{Evolution of the objective function for the different algorithms for matrix pair $(E,A)$ for the perturbed mass-spring-damper system. \label{MSD101ton}}
\end{figure*}

The corresponding mass vector $m$ containing the point masses,
spring vector $k$ containing the spring constants, and damping vector $c$ containing the
damping constants are all equal to the vector
\[
m=c=k=\mat{cccccccccc}1 &2&3&4&5&6&7&8&9&10 \rix.
\]
Figure~\ref{MSD101ton} shows the evolution of the objective function~\eqref{reformoptprob}. We used the initialization from~\eqref{msdmat}, i.~e. we used the true $J$, $R$ and $Q$ from the known stable system to see whether an approximation better than the original stable pencil can be found. (Note that if we use our standard initialization, see Section~\ref{subsecalgo}, FGM converges to another local minimum with higher approximation error $32.70$.)

We observe that FGM converges much faster than GM, as in the two other examples,
while it is able to generate a better approximation than the FGM from~\cite{GilS16} applied on the nearest stable matrix problem $(I,E^{-1} A)$,
as already noted in Section~\ref{case1}.
We also observe that the algorithms are able, rather surprisingly, to significantly reduce the value of the objective compared to the initialization (which is the original stable pencil),
from $21.97$
to $6.53$ for FGM (DH matrix),
to $12.70$ for GM (DH pair), and
to $4.09$ for FGM (DH pair). Actually, if we run FGM for $5$ more seconds, then it converges to a solution with a value $3.81$.
Figure~\ref{MSD101toneig} shows the location of the eigenvalues of the pencils $(E,A)$ generated by the FGM using a DH matrix and a DH pair, as well as those of the original and perturbed pencils.
\begin{figure*}[h!]
\centering
\includegraphics[width=\textwidth]{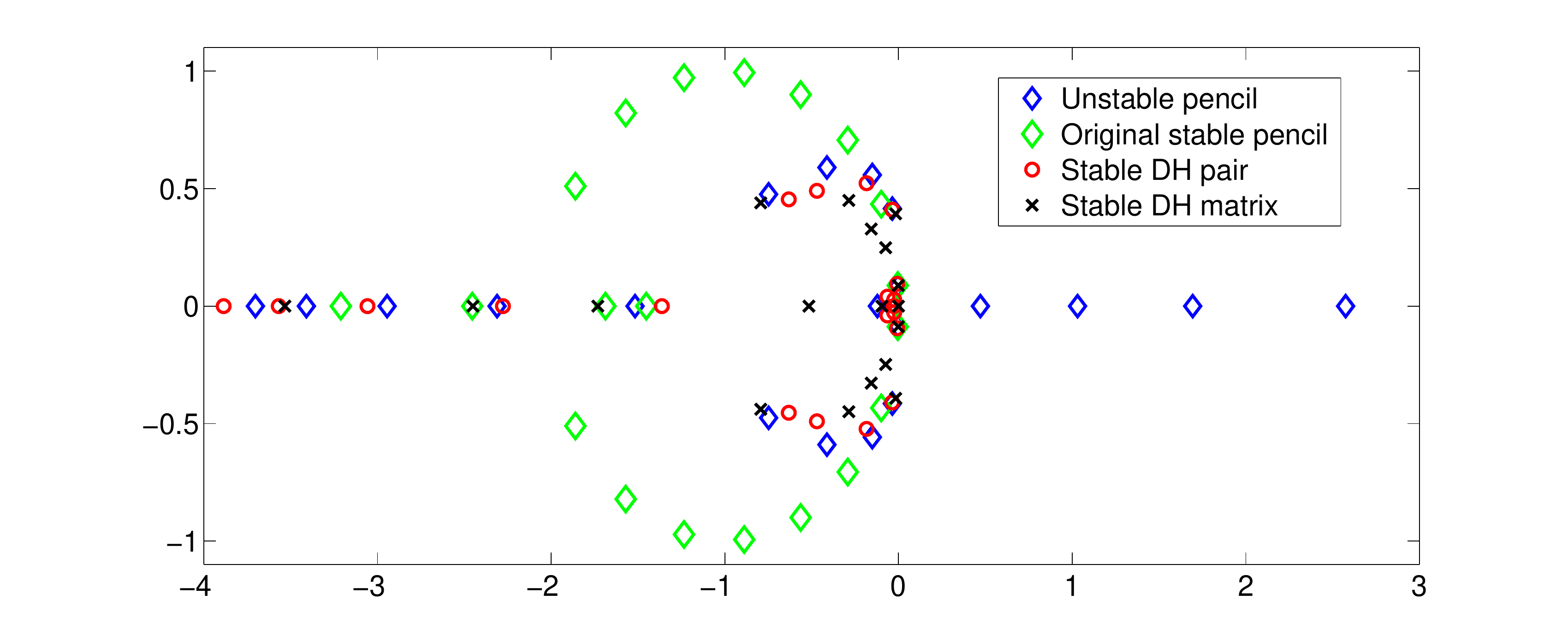}
\caption{Eigenvalues of various pencils corresponding to the mass-spring-damper system. \label{MSD101toneig}}
\end{figure*}

\section*{ Conclusion }

We have considered the problem of computing the nearest stable matrix pair to a given unstable one.
We have presented a new characterization of the set of stable matrix pairs using DH matrix pairs, defined as $(E,(J-R)Q)$ where
$J=-J^T$, $R \succeq 0$, and  $Q^TE \succeq 0$ with $Q$ invertible. This has allowed us to reformulate the
nearest stable matrix pair problem in a more convenient way for which we develop an efficient fast
gradient method.

Possible future work is related to develop more sophisticated ways to solve~\eqref{reformoptprob} and to apply
globalization approaches to obtain better solutions.

\bibliographystyle{siam}
\bibliography{GilMS}

\newpage

\appendix

\section{ Gradient of $f$ with respect to $Q$ } \label{appgradQ}

In this section, we explain how to compute the gradient of
\[
f(J,R,Q,H) = \frac{1}{2} {\|A - (J-R) Q\|}_F^2 + \frac{1}{2} {\|E - Q^{-T}H\|}_F^2,
\]
with respect to $Q$, denoted by $\nabla_Q f(Q)$.
For the first term, we have
\[
\nabla_Q \frac{1}{2} {\|A - (J-R) Q\|}_F^2 = -(J-R)^T (A - (J-R) Q).
\]
For the second term, we need to compute $\nabla_Q {\|E- Q^{-T}H\|}_F^2$.
Using the basic  rules for matrix differentiation, see e.~g.~\cite{MagN99, olsen2012efficient}, 
\begin{enumerate}
\item[(i)] Linear:
\[
\nabla_X \tr \big( WX \big)
=
W^T,\quad
\nabla_X \tr \big( WX^T \big)
=
A.
\]
\item[(ii)] Product:
\[
\nabla_X \tr \big( G(X) F(X) \big)
=
\nabla_X \tr \big( F(Y) G(X)  + F(X) G(Y) \big) |_{Y=X} .
\]

\item[(iii)] Inverse:
\[
\nabla_X \tr \big( W F^{-1}(X) \big)
=
- \nabla_X \tr \big( F^{-1}(Y) W F^{-1}(Y) F(X)  \big) |_{Y=X},
\]
\end{enumerate}
we obtain (considering the transpose matrix $E^T - HQ^{-1}$),
\[
\nabla_Q  {\|E^T - HQ^{-1}\|}_F^2
 =
\nabla_Q \tr \big( (E^T - HQ^{-1})^T (E^T - HQ^{-1}) \big)
 =
\nabla_Q \tr ( Q^{-T} H^T H Q^{-1}  ) - 2 \tr ( E H Q^{-1} ) .
\]
Using (ii) with $F(Q) = H Q^{-1}$ and $G(Q) = F(Q)^T$, we get for the first term
\begin{align*}
\nabla_Q \tr ( G(Q) F(Q)  )
& =
\nabla_Q \tr ( F(Y) G(Q)  + F(Q) G(Y)  )  |_{Y=Q} \\
& =
\nabla_Q \tr ( F(Y) F(Q)^T  + F(Q) F(Y)^T  )  |_{Y=Q} \\
& =
2 \nabla_Q \tr ( F(Q) F(Y)^T  )  |_{Y=Q} \\
& =
2 \nabla_Q \tr ( H Q^{-1}  Y^{-T} H^T  )  |_{Y=Q} .
\end{align*}
Using (ii) again, now with $G(Q) = H Q^{-1}$ and $F(Q) = Y^{-T} H^T$ which is a constant, we obtain
\[
\nabla_Q \tr ( H Q^{-1} Y^{-T} H^T  )
=
\nabla_Q \tr ( Y^{-T} H^T H Q^{-1}  + Y^{-T} H^T H Y^{-1}  )
=
\nabla_Q \tr ( (Y^{-T} H^T H) Q^{-1} )  ,
\]
using (iii),  $W = Y^{-T} H^T H$, and (i),
\[
\nabla_Q \tr ( (Y^{-T} H^T H) Q^{-1} ) |_{Y=Q}
= -( Q^{-1} Q^{-T} H^T H Q^{-1}  )^T
= - Q^{-T} H^T H Q^{-1} Q^{-T}.
\]
For the second term we use (iii) and obtain
\[
\nabla_Q \tr ( E H Q^{-1} )
=
-( Q^{-1} E H Q^{-1} )^T
=
- Q^{-T} H^T E^T Q^{-T}.
\]

\end{document}